\documentclass[12pt]{elsarticle}

\makeatletter
\def\@preprinttext{}  
\def\@preprint{\@preprinttext}  
\makeatother

\usepackage{amssymb,hyperref}
\usepackage{amsmath}
\usepackage{multirow}
\usepackage{algorithm}
\usepackage{bm}
\usepackage{algorithmic}
\usepackage{array} 
\usepackage{color}
\usepackage{caption} 
\usepackage{graphicx}   
\usepackage{subfigure}

\journal{}

\begin{document}

\begin{frontmatter}



\title{Meta-Learning the Optimal Mixture of Strategies for Online Portfolio Selection}


\author{
    Jiayu Shen,
    Jia liu\footnote{Corresponding author. Email: jialiu@xjtu.edu.cn}, 
    Zhiping Chen
} 


\affiliation{organization={School of Mathematics and Statistics, Xi'an Jiaotong University},
            addressline={Xianning West Road 28}, 
            city={Xi'an},
            postcode={710049}, 
            state={Shaanxi},
            country={P. R. China}}

\begin{abstract}
This paper presents an innovative online portfolio selection model, situated within a meta-learning framework, that leverages a mixture policies strategy. The core idea is to simulate a fund that employs multiple fund managers, each skilled in handling different market environments, and dynamically allocate our funding to these fund managers for investment. To address the non-stationary nature of financial markets, we divide the long-term process into multiple short-term processes to adapt to changing environments. We use a clustering method to identify a set of historically high-performing policies, characterized by low similarity, as candidate policies. Additionally, we employ a meta-learning method to search for initial parameters that can quickly adapt to upcoming target investment tasks, effectively providing a set of well-suited initial strategies. Subsequently, we update the initial parameters using the target tasks and determine the optimal mixture weights for these candidate policies. Empirical tests show that our algorithm excels in terms of training time and data requirements, making it particularly suitable for high-frequency algorithmic trading. To validate the effectiveness of our method, we conduct numerical tests on cross-training datasets, demonstrating its excellent transferability and robustness.
\end{abstract}



\date{}
\end{frontmatter}



\section{Introduction}
\label{sec1}
Long-term portfolio management is a fundamental problem in finance, involving the periodic reallocation of assets. 
Online portfolio selection (OLPS) is 
an important and challenging problem in 
long-term portfolio management, which aims to maximize the expected returns while controlling risks through reallocating portfolio weights. 
   
According to the survey \cite{Survey}, existing OLPS algorithms can be broadly categorized into the following five classes: benchmark algorithms, follow-the-winner (FTW) algorithms, follow-the-loser (FTL) algorithms, pattern-matching algorithms, and mixture policies.
Traditional online portfolio algorithms mostly rely on artificially designed financial features. For example, FTL algorithms assume that the market follows the mean reversion principle. Building on this idea, Passive Aggressive Mean Reversion (PAMR) introduces an innovative design for its loss function \cite{PAMR}: if the expected return of the price in the previous period exceeds a specified threshold, the loss gradually increases; conversely, if the expected return is below or equal to the threshold, the loss is zero. In contrast, the FTW algorithm assumes that stocks currently performing well will continue to perform well. For instance, the Exponentiated Gradient (EG) strategy seeks to track the best-performing stock from the previous period while maintaining the portfolio as close as possible to its previous allocation \cite{EG}. Due to these assumptions, {\it most traditional investment strategies perform well only in specific {datasets} and can hardly adapt to the change of market or the market regime}. Unlike the other four algorithms, mixture policies 
do not focus solely on a single strategy but maintain a pool of candidate policies and allocate funds among them. This paper adopts this idea by maintaining a pool of candidate policies and determining the allocation ratios of these policies based on current market conditions.

More recently, deep learning has been developed to address the portfolio management problem\cite{MA2021113973,martinez2024portfolio,2021End,niu2022metatrader}. {It provides new insights into online portfolio management, avoids various assumptions about the market, and achieves objective profitability.}
{\it However, deep learning encounters two challenges in online portfolio selection problems: long training time and high sample complexity}. Financial markets are highly dynamic, requiring investors to make timely decisions, particularly real time decisions for online portfolio selection problems. Additionally, due to the high instability of financial data sequences and quick variation of the corresponding distribution over time, it is difficult to make a precise prediction. Meanwhile, the historical data contains a vast amount of irrelevant information, with only recent data being of significance. Training and maintaining a end-to-end neural network demands extensive datasets and high time costs, which are challenging to satisfy the real time requirement in online portfolio selection problems.

Meta-learning, 
an emerging machine learning method, focuses on rapidly learning, generalizing, and adapting across different tasks \cite{vilalta2002perspective,hospedales2021meta,li2019meta,shu2019meta}.
The characteristics of meta-learning make it naturally suitable
for online portfolio management problems. By learning patterns and rules from historical data, meta-learning enables models to better adapt to unknown environments and new investment scenarios. Its core strength lies in rapidly learning and adapting to new market features and conditions, thereby enhancing the robustness and generalization of investment portfolios.
In recent years, an increasing number of papers have applied meta-learning to time series forecasting and stock trading problems \cite{woo2022deeptime,Ma2023,2020Meta,10.1145/3449639.3459386}. {\it However, due to differences in data structures across various markets, these methods lack the ability to transfer effectively between markets}. This limitation results in meta-learning models being able to learn from only a single market, which restricts the model's generalization ability.

Motivated by the aforementioned insights, {\it this paper integrates meta-learning approach with a mixture policies learning framework, resulting the Meta-LMPS-online model, to make online portfolio selection.} We decompose the long-term investment process into multiple short-term time segments, generating several small-sample tasks to address the challenge of evolving data distributions in financial markets over time. During training, we employ clustering techniques to identify historically high-performing and diverse policies as candidate policies, and utilize meta-learning approach to derive initial parameters capable of rapidly adapting to new tasks. This approach enables the model to swiftly and dynamically adjust investment strategies in response to new tasks, while demonstrating robust transferability and generalization across various market conditions.

Our main contributions are summarized in three points.
\begin{itemize}
\item Unlike traditional deep learning or meta-learning frameworks that directly learn the portfolio, we adopt a mixture policy learning framework, which focuses on learning the mixture weights of multiple online portfolio policies independent of the number of stocks. This approach effectively addresses the transferability and generalization challenges across diverse markets with varying stock pool sizes and datasets.
\item We address the challenge of non-stationarity in financial markets by decomposing the long-term investment process into multiple short-term tasks. 
\item Numerical results valid the superiority of the proposed algorithm over traditional OLPS strategies, when transferring between different time periods or markets.
\end{itemize}

In Section 2, we introduce the background of online portfolio selection and related work. In Section 3, we introduce the model-agnostic meta-learning for online portfolio selection. In Section 4, we carry out a series of numerical tests. Section 5 concludes.
Table \ref{Symbol Glossary} introduce the mathematical notations in the paper.

\begin{table}[h!]
\centering
\caption{
 Notations
 }
	\label{Symbol Glossary}{%
    \resizebox{\textwidth}{!}{ 
\begin{tabular}{cc}
		\hline
 Symbol & Definition or Explanation\\ \hline
			
        $N$ & Number of assets\\
        $T$ & Number of trading periods\\
        $\boldsymbol{p}_t=(p_{t,1},\dots,p_{t,N})^\top$ & Closing price at the end of period $t$\\
        $\boldsymbol{x}_t=\boldsymbol{p}_t/\boldsymbol{p}_{t-1}$ & Price relative vector of period $t$\\
        $\boldsymbol{x}_{[t_1,t_2]}$ & Price changes over period $t_1$ to period $t_2$\\
        $\boldsymbol{b}_t=(b_{t,1},\dots,b_{t,N})^\top\in  \Delta_N^+$ & Portfolio vector/Capital allocation of assets at the start of period $t$\\
        $\boldsymbol{b}_{[t_1,t_2]}$ & Portfolio from period $t_1$ to period $t_2$ \\
        $\widetilde{\boldsymbol{b}}_t=\frac{\boldsymbol{b}_t \odot \boldsymbol{x}_t}{\boldsymbol{b}_t^\top \boldsymbol{x}_t}$ & Capital allocation of assets at the end of period $t$\\
        $S_t$ & Cumulative wealth at the end of period($S_0$ denote the initial wealth)\\
        $\mu_{t-1}$ & Transaction remainder factor(TRF)\\
        $\delta$ & Unit transaction cost for buying/selling assets\\
            $M$ & Number of candidate policies\\
        $\boldsymbol{r}_t$ & Return rate of candidate policies in period $t$\\
        $w$ & Historical window size within each sub-data pair in the task\\
        $\boldsymbol{\omega}_t=(\omega_t^1,\dots,\omega_t^M)\in \Delta_M^+$ &Mixture weight in period $t$\\
        $f_\theta$ & Neural network model\\
        $\mathcal{T}_l$ & Sub-task $l$\\
        $\mathcal{T}_l^{sup}$ & Support Set of $\mathcal{T}_l$\\
        $\mathcal{T}_l^{que}$ & Query Set of $\mathcal{T}_l$\\
        $K$ & Number samples in $\mathcal{T}_l^{sup}$\\
        $Q$ & Number samples in $\mathcal{T}_l^{que}$\\
        $H$ & Number of batch training tasks\\
        $\alpha$ & Inner loop learning rate\\
        $\beta$ & Outer loop learning rate\\
        \hline
    \end{tabular}%
}
}
\end{table}

\section{Online portfolio selection and related work}

\textbf{Problem setting.}
We consider the problem that an investor plans to allocate funds across a finite number of assets over a series of time periods. Specifically, we consider $ N \geq 2 $ assets and $ T \geq 1 $ periods.
Let $ \boldsymbol{p}_t = (p_{t,1}, p_{t,2}, \dots, p_{t,N}) \in \boldsymbol{R}_+^N $ denote the closing prices of the $ N $ assets at the end of period $ t $. These closing prices serve as the opening prices for the assets at the start of period $ t+1 $.
In period $t$, we consider the price relative vector $ \boldsymbol{x}_t = (x_{t,1}, x_{t,2}, \dots, x_{t,N}) \in \boldsymbol{R}_+^N $, where $x_{t,i} = \frac{p_{t,i}}{p_{t-1,i}}$ representing the growth factor of asset $ i $ during period $ t $. The market price changes from period $ t_1 $ through period $ t_2 $ are represented by $ \boldsymbol{x}_{[t_1,t_2]} = \{ \boldsymbol{x}_{t_1}, \dots, \boldsymbol{x}_{t_2} \} $.
At the beginning of period $t$, the investment is allocated according to the portfolio vector $ \boldsymbol{b}_t = (b_{t,1}, b_{t,2}, \dots, b_{t,N})^\top $, where $ b_{t,i} $ denotes the proportion of capital allocated to asset $ i $. Assuming a self-financing investment process without allowing short-selling positions, the portfolio vector adheres to the constraints $ \boldsymbol{b}_t \in \Delta_N^+ $, where $ \Delta_N^+ = \{ \boldsymbol{b} \mid \boldsymbol{b} \geq \boldsymbol{0}, \, \boldsymbol{b}^\top \boldsymbol{1} = 1 \} $. The portfolio strategy from period $ t_1 $ through period $ t_2 $ is denoted by $ \boldsymbol{b}_{[t_1,t_2]} = \{ \boldsymbol{b}_{t_1}, \dots, \boldsymbol{b}_{t_2} \} $.

Let $ S_t $ denote the cumulative wealth of the investor at the beginning of period $ t $. After the market fluctuations during period $ t $, the change in the investor's wealth is $ S_t \boldsymbol{b}_t^\top\boldsymbol{x}_t $. The portfolio changes are represented by $\widetilde{\boldsymbol{b}}_t=(\widetilde{b}_{t,1},\dots,\widetilde{b}_{t,N})^\top$, where $\widetilde{b}_{t,i}=\frac{b_{t,i}x_{t,i}}{\boldsymbol{b}_{t}^\top\boldsymbol{x}_{t}}$. At the beginning of period $ t+1 $, the investment strategy is adjusted to $\boldsymbol{b}_{t+1}$, resulting in transaction costs. Let $\mu_t$ denote the  transaction remainder factor \cite{magill1976portfolio,davis1990portfolio,li2018transaction}, which represents the remaining proportion of the capital after the deduction of transaction cost at the beginning of period $t$. Consequently, the change in the investor's wealth at this time is $\mu_t S_t \boldsymbol{b}_t^\top\boldsymbol{x}_t$. Therefore, after $ T $ periods, the investor's cumulative wealth is computed as
\begin{align}
\label{culmulative wealth}
S_T(\bm{b}_{[1,T]})=S_0\prod_{t=1}^T\mu_{t-1}\boldsymbol{b}_t^\top\boldsymbol{x}_t.
\end{align}
Assuming the transaction cost rate is $\delta$, we approximate Equation \eqref{culmulative wealth} as follows:
{\small
\begin{align*}
    S_T(\boldsymbol{b}_{[1,T]}) \approx & S_0 \prod_{t=1}^T \Bigg[ \sum_{i=1}^N b_{t,i}x_{t,i}\left(1-\delta \sum_{i=1}^N \left| b_{t,i} - \widetilde{b}_{t-1,i} \right|\right)\Bigg].
\end{align*}}
This approximation technique is proposed by \citet{GUO2023}.

\textbf{Traditional OLPS methods.} Table \ref{tab:policy}
 lists classical online portfolio selection policies, including benchmarks, FTW, FTL, and pattern-matching categories. We use these algorithms to build our \textbf{candidate policy pool}. For mixture policies, \citet{FU} allocates assets evenly a mong experts who operate autonomously, then aggregates their wealth. \citet{FLH} maintains a limited set of experts, dynamically adds or removes them based on performance, and adjusts weights among the active experts. 
 \citet{lin2024online} and \citet{yang2024aggregating} apply the online gradient update algorithm or exponentially weighted average algorithm to integrate a pool of expert strategies. 
 These approach smooths overall performance by balancing various policies and enhances adaptability by combining universal strategies with heuristic algorithms.

\begin{table}[htp]
    \caption{Classical Online Portfolio Strategies}
    \scalebox{0.9}[0.9]{\label{tab:policy}
    \begin{tabular}{|>{\centering\arraybackslash}m{2.4cm}|>{\centering\arraybackslash}m{1.3cm}|>{\centering\arraybackslash}m{2.2cm}|>{\centering\arraybackslash}m{8.3cm}|}
    \hline
    Classifications & Number & Strategies & Updating rule of $\boldsymbol{b}_{t+1}$ \\ \hline
    \multirow{2}{*}{Benchmarks}  
    & 1 & BAH\cite{Survey} & $\frac{\boldsymbol{b}_t\bigotimes \boldsymbol{x}_t}{\boldsymbol{b}_t\cdot \boldsymbol{x}_t}$ \\ \cline{2-4} 
    & 2 & CRP\cite{CRP} & $\boldsymbol{b}_t$\\ \hline
    \multirow{4}{*}{FTW}  
    & 3 & UP\cite{UP} & $\frac{\int_{\Delta_N^+} \boldsymbol{b}S_t(\boldsymbol{b})du(\boldsymbol{b})}{\int_{\Delta_N^+} S_t(\boldsymbol{b})du(\boldsymbol{b})}$ \\ \cline{2-4} 
    & 4 & Best-so-far &  $\operatorname*{argmax}_{\boldsymbol{b} \in \Delta_N^+} \sum_{i=0}^{w-1} \frac{1}{w}\text{log} (\boldsymbol{b}\cdot\boldsymbol{x}_{t-w})$ \\ \cline{2-4}
    & 5 & EG\cite{EG} & $\operatorname*{argmax}_{\boldsymbol{b} \in \Delta_N^+} \eta \text{log}(\boldsymbol{b}\cdot \boldsymbol{x}_t)-R \left(\boldsymbol{b},\boldsymbol{b}_t\right)$ \\ \cline{2-4} 
    & 6 & ONS\cite{ONS} & $\operatorname*{argmax}_{\boldsymbol{b} \in \Delta_N^+} \sum_{i=1}^t  \text{log} \left(\boldsymbol{b}\cdot \boldsymbol{x}_i\right) - \lambda R \left(\boldsymbol{b}\right)$ \\ \hline
    \multirow{6}{*}{FTL}
    & 7 & ANTICOR\cite{ANTICOR} & If assets $i$ increases more than asset $j$ and they are positively correlated, Anticor transfers from $i$ to $j$ by $M_{cor}(i,j)-\text{min}\{0,M_{cor}(i,i)\}-\text{min}\{0,M_{cor}(j,j)\}$ (where $M_{cor}$ is the cross correlation value) \\ \cline{2-4} 
    & 8 & PAMR\cite{PAMR} & $\rm{arg} \rm{min}_{\boldsymbol{b} \in \Delta_N^+} R \left(\boldsymbol{b},\boldsymbol{b}_t\right)$ s.t.  $\boldsymbol{b} \cdot \boldsymbol{x}_{t}\leq \epsilon $ \\ \cline{2-4} 
    & 9 & CWMR\cite{CWMR} & $(\boldsymbol{\mu}_{t+1},\boldsymbol{\sum}_{t+1})=\operatorname*{argmin}_{\boldsymbol{\mu},\boldsymbol{\sum}}D_{KL}(\mathcal{N}(\boldsymbol{\mu},\boldsymbol{\sum})||\mathcal{N}(\boldsymbol{\mu}_t,\boldsymbol{\sum}_t)$ s.t. $Pr(\boldsymbol{\mu} \cdot \boldsymbol{x}_t\leq \epsilon)\geq \theta$ $\boldsymbol{\mu}\in \Delta_N, \boldsymbol{b}_{t+1}\in \mathcal{N}(\boldsymbol{\mu}_{t+1},\boldsymbol{\sum}_{t+1})$ \\ \cline{2-4} 
    & 10 & OLMAR\cite{OLMAR} & $\operatorname*{argmin}_{\boldsymbol{b} \in \Delta_N^+} R \left(\boldsymbol{b},\boldsymbol{b}_t\right)$ s.t. $\boldsymbol{b} \cdot \hat{\boldsymbol{x}}_{t+1}\geq \epsilon$ $\hat{\boldsymbol{x}}_{t+1}(w)=\frac{1}{w}\left( 1+\frac{1}{\boldsymbol{x}_t}+\dots+\frac{1}{\prod_{i=0}^{w-2}\boldsymbol{x}_{t-i}}\right)$ \\ \cline{2-4} 
    & 11 & WMAMR\cite{WMAMR} & $\operatorname*{argmin}_{b\in \nabla_N^+} \frac{1}{2} ||\boldsymbol{b}-\boldsymbol{b}_t||$ s.t. $\boldsymbol{b} \cdot \widetilde{\boldsymbol{x}}_{t,w}\leq \epsilon$,
    $\widetilde{\boldsymbol{x}}_{t,w}=\frac{1}{w}\left( 1+\frac{1}{\boldsymbol{x}_t}+\dots+\frac{1}{\prod_{i=0}^{w-2}\boldsymbol{x}_{t-i}}\right)$ \\ 
    \cline{2-4}
    &12 &RMR\cite{RMR} & $\operatorname*{argmin}_{\boldsymbol{b} \in \Delta_N^+} \boldsymbol{b} \cdot \widehat{\boldsymbol{x}}_{t}+\lambda R \left(\boldsymbol{b},\boldsymbol{b}_t\right)$ $\widehat{\boldsymbol{x}}_{t}=\frac{\boldsymbol{u}_{t+1}}{\boldsymbol{p}_{t}}, \boldsymbol{u}_{t+1}=\text{min}_u\Sigma_{i=0}^{\omega-1}||\boldsymbol{p}_{t-i}-\boldsymbol{u}||$ \\
    \hline
    \multirow{2}{*}{\begin{tabular}[c]{@{}c@{}}Pattern\\ Matching\end{tabular}}  
    & 13 & $\text{B}^{\text{NN}}$\cite{BNN} & \qquad \quad$\operatorname*{argmax}_{\boldsymbol{b} \in \Delta_N^+} \sum_{i\in \emph{C}_t} \frac{1}{|\emph{C}_t|}\text{log} \boldsymbol{b}\cdot\boldsymbol{x_i}$\quad \qquad ($\emph{C}_t$ is the set of similar windows; pattern matching methods use different criteria) \\ \cline{2-4}
    & 14 & CORN\cite{CORN} & $\operatorname*{argmax}_{\boldsymbol{b} \in \Delta_N^+} \sum_{i\in \emph{C}_t} \frac{1}{|\emph{C}_t|}\text{log} \boldsymbol{b}\cdot\boldsymbol{x_i}$ \\  
    \hline
    \end{tabular}
    }
\end{table}

\textbf{Deep learning in portfolio optimization.} 
The application of deep learning in OLPS can be divided into two categories. The first category involves using deep learning to predict asset returns and subsequently constructing an optimization model to formulate portfolios based on the predictions. The second category adopts an end-to-end approach, where deep neural networks are used to directly generate portfolios.
\citet{2017Threshold, MA2021113973, martinez2024portfolio} use recurrent neural networks or convolutional networks to predict stock prices and employ optimization methods to determine investment strategies. \citet{2021End, 9403044, zhang2021universalendtoendapproachportfolio} combine prediction and optimization tasks into an end-to-end neural network, considering various metrics including investment returns and volatility as input features for the network. \citet{niu2022metatrader} explored the application of mixture policy in deep learning, integrating expert policies and environmental interactions to adapt to dynamic markets, providing a new perspective for online portfolio optimization.


\textbf{Meta learning in time series forecasting and portfolio management.} \citet{DBLP} and \citet{woo2022deeptime} applied meta-learning to historical value models and time index models respectively to tackle the problem of non-stationary forecasting, which verified the improvement of model generalization by meta-learning. \citet{2020Meta} applied a quantum meta-learning algorithm to the field of stock trading, proposing a quantum neural network model. They demonstrated through simulations on quantum computers that it can achieve efficient learning and adapt trading strategies even with a significant reduction in the number of parameters. \citet{10.1145/3449639.3459386} applied meta reinforcement learning to the online portfolio optimization problem, by breaking down long-term investment goals into a series of short-term tasks
and utilizing convolutional networks to match and evaluate market data patterns.
However, these methods have not yet fully leveraged the potential of meta-learning in the field of online portfolio management, and they lack the ability to transfer and generalize across multiple markets.

\section{Model-agnostic
meta-learning for online portfolio selection}
\begin{figure}[ht]
    \centering
    \includegraphics[width=0.7\textwidth]{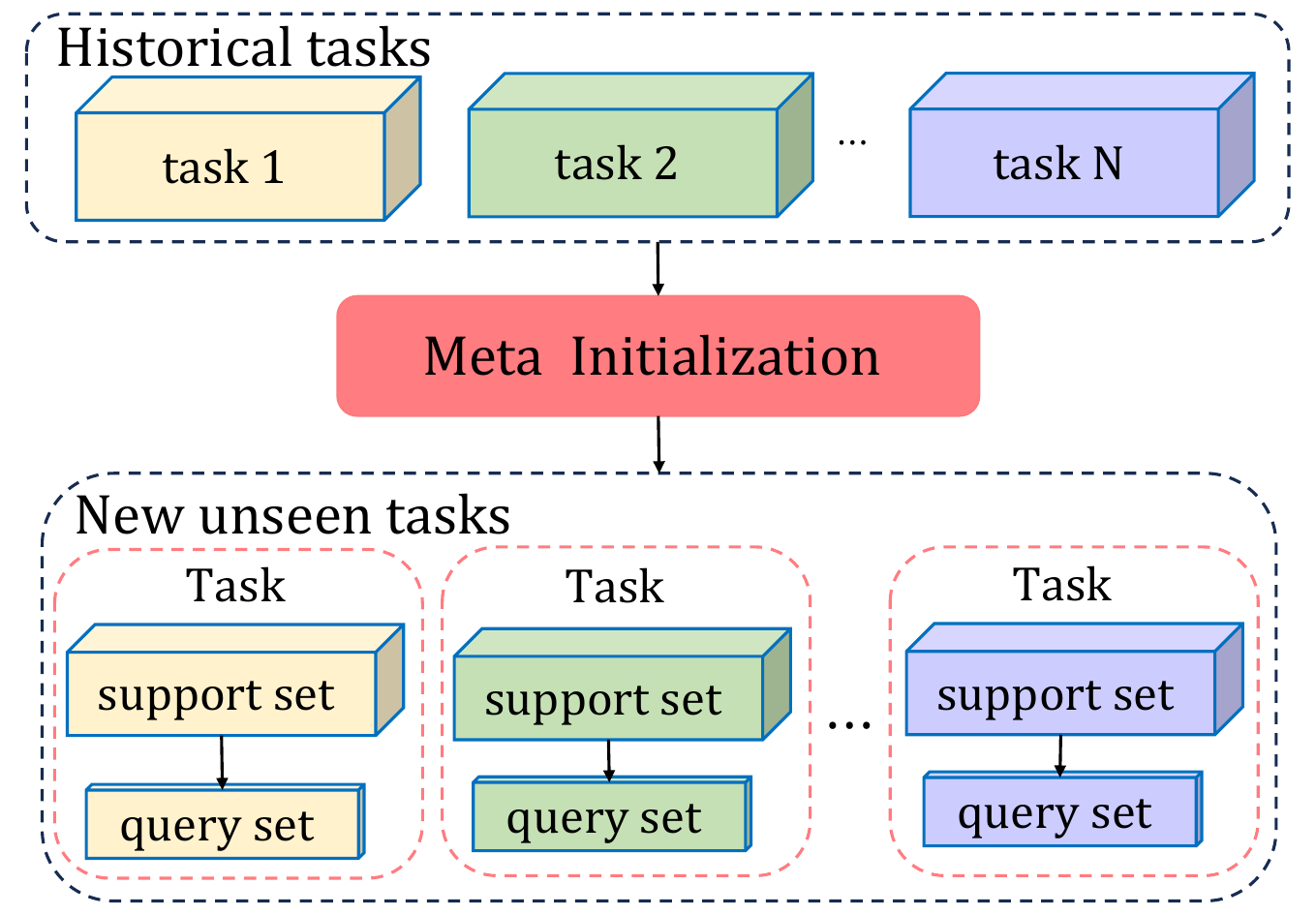}
    \caption{Meta learning}
    \label{fig:Meta}
\end{figure}

The meta-learning framework employed in this paper is model-agnostic
meta-learning (MAML) proposed by \citet{MAML}. The core concept involves learning an initial set of parameters from historical tasks that can quickly adapt to new tasks. The background introduction of MAML is provided in Appendix A. The process of learning from historical tasks is known as meta-training and is performed offline, whereas the adaptation to new tasks is known as meta-testing and is performed online (ref. Figure \ref{fig:Meta}).
\subsection{Single Task in Online Portfolio Selection} 





{\textbf{Portfolio selection task.}
Considering the high instability of financial data sequences and quick variation of the corresponding distribution over time, we consider the investment process in a short time period in one single market as a single task. Then the investment tasks over several segments in one single market or in different markets can be viewed as a meta task.

Suppose in a given market, we utilize the historical relative price data of stocks over the previous $L$ periods, $\boldsymbol{x}_{[1,L]}=\{\boldsymbol{x}_1,\boldsymbol{x}_2,\dots,\boldsymbol{x}_L\}$, to construct the investment portfolio vector for the subsequent $Q$ periods. From the policy pool (Table  \ref{tab:policy}), we select $M$ candidate policies that demonstrate strong performance and diverse decision-making under varying market conditions. Specifically, we apply k-means clustering to the return sequences of policies from the policy pool over the first $L$ periods. The policy with the highest cumulative return in each cluster is selected as a candidate.

Each candidate policy generates an investment portfolio vector for every period based on the historical relative price sequence. Assuming that the investment portfolios of $M$ candidate policies at the start of period $t$ are $\boldsymbol{b}_t^1,\dots,\boldsymbol{b}_t^M$, and their return rate at period $t$ are $\boldsymbol{r}_t=(r_t^1,\dots,r_t^M).$ 
Different from traditional deep learning or meta-learning
frameworks, we learn the mixture weights of the candidate portfolio policies 
$\boldsymbol{\omega}_t$ rather than the portfolio $\boldsymbol{b}_t$ directly. We set the historical window size as $w$.
 A single task of online portfolio selection is to determine the optimal mixture weights $\boldsymbol{\omega}_{t}=({\omega}^1_{t},\dots,{\omega}^M_{t})$ among $M$ candidate policies, for the following period based on the historical return data $[\boldsymbol{r}_{t-w},\dots,\boldsymbol{r}_{t-1}]$, $t=L+1,\ldots,L+Q$. 

We then denote the portfolio selection model in a task as 
\begin{align}
\label{eq:network_mapping}
f_\theta: [\boldsymbol{r}_{t-w},\dots,\boldsymbol{r}_{t-1}]\rightarrow \boldsymbol{\omega}_{t},\  t=L+1,\dots,L+Q,
\end{align}
where $\theta$ is the parameter of the training model.
After obtaining the mixture weights through the neural network, we can ultimately derive the investment portfolio vector as
\begin{align*}
\boldsymbol{b}_{t}=\sum_{j=1}^M \omega_{t}^j \boldsymbol{b}_{t}^j,\quad t=L+1,\dots,L+Q.
\end{align*}
}

\textbf{Neural Network.}
We then design the neural network. As shown in Equation \eqref{eq:network_mapping}, the input of the network consists of the historical return rates of the candidate policies $[\boldsymbol{r}_{t-w},\dots,\boldsymbol{r}_{t-1}]$ over the previous $w$ periods. The output of the network is the mixture weights $\boldsymbol{\omega}_{t}$, and its parameters are denoted as $\theta$.

\begin{figure}[h!]
    \centering
    \includegraphics[width=\textwidth]{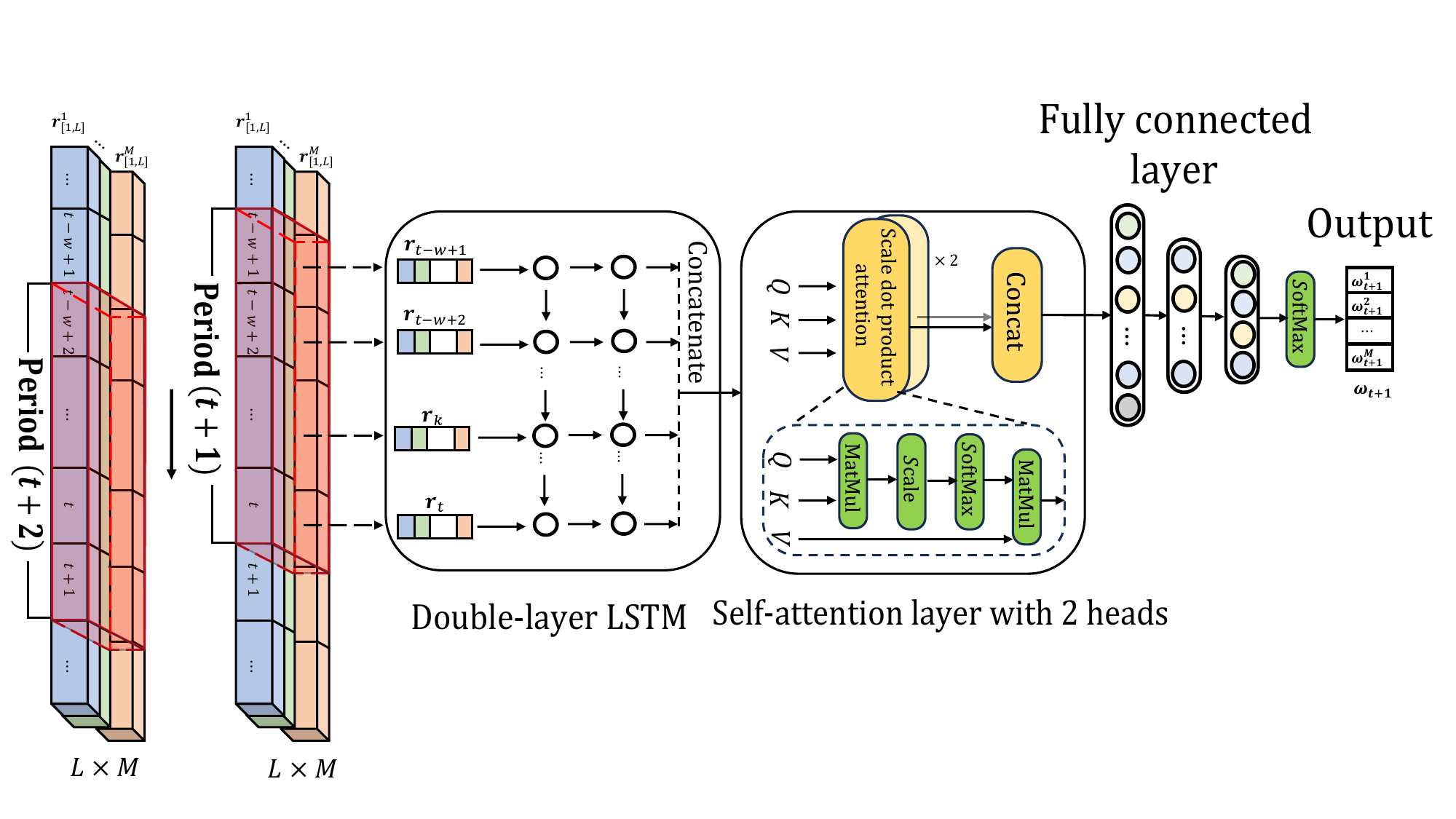}
    \caption{Neural network structure}
    \label{fig:Network structure}
\end{figure}

The network architecture is shown in Figure \ref{fig:Network structure} and comprises three components: a two-layer LSTM, an attention mechanism, and fully connected layers \cite{cipiloglu2022portfolio,rather2021lstm,alzaman2024deep}. The two-layer LSTM is employed to extract features from the time series data and capture sequential dependencies. The self-attention mechanism is utilized to discern correlations among the vectors. The fully connected layers are responsible for mapping the extracted features to the final output space, with the Softmax layer providing the probability distribution over the selected candidate policies.

\textbf{Data and labels.}
To train the neural network, we need $L-w$ pairs of labeled data:
$$\{[\boldsymbol{r}_{t-w},\dots,\boldsymbol{r}_{t-1}];\hat{\boldsymbol{\omega}}_{t}\}_{t=w+1}^{L}.$$
Here, 
$\hat{\boldsymbol{\omega}}_{t}=(\hat{\omega}^1_{t},\dots,\hat{\omega}^M_{t})$ represents a label vector composed of the relative weights of $M$ candidate policies. 
The weights are determined by the actual returns of each policy in the period $t$. 
The specific calculation for $\hat{\boldsymbol{\omega}}_{t}$ is as follow:
\begin{align*}
    \hat{\omega}^j_{t} &= \frac{r_{t}^j - \text{min}\{r_{t}^1, \dots, r_{t}^M\}}{\sum_{j=1}^M r_{t}^j - M \times \text{min}\{r_{t}^1, \dots, r_{t}^M\}}, \\
    & \qquad \qquad \quad j=1, \dots, M, \quad t=w+1, \dots, L.\notag
\end{align*}


\textbf{Loss function.}
Given the label $\hat{\boldsymbol{\omega}}_{t}$ at period $t$, 
we define the loss function as:
\begin{align}
    \mathcal{L}
(\boldsymbol{\omega}_{t},\hat{\boldsymbol{\omega}}_{t})&=\Vert \boldsymbol{\omega}_{t}-\hat{\boldsymbol{\omega}}_{t}\Vert_2^2+\eta \Vert \boldsymbol{\omega}_{t}\Vert_2^2 \notag\\
\qquad \qquad+&\lambda \delta \Vert \sum_{j=1}^M (\omega_{t}^j \boldsymbol{b}_{t}^j-\omega_{t-1}^j \widetilde{\boldsymbol{b}}_{t-1}^j)\Vert_1 .
\end{align}

The loss function consists of three terms. The first term is used to control the gap between the predicted values and the actual values, the second term is for encouraging a diversified allocation of weights among each strategy, and the third term is used to control transaction costs. Here, $\eta$ and $\lambda$ are both positive constants.


\subsection{
Meta-Learning Framework}
We then consider the meta task with many sub-investment tasks over several time segments in one single market or in several different markets.

\begin{figure}[ht]
    \centering
    \includegraphics[width=0.9\textwidth]{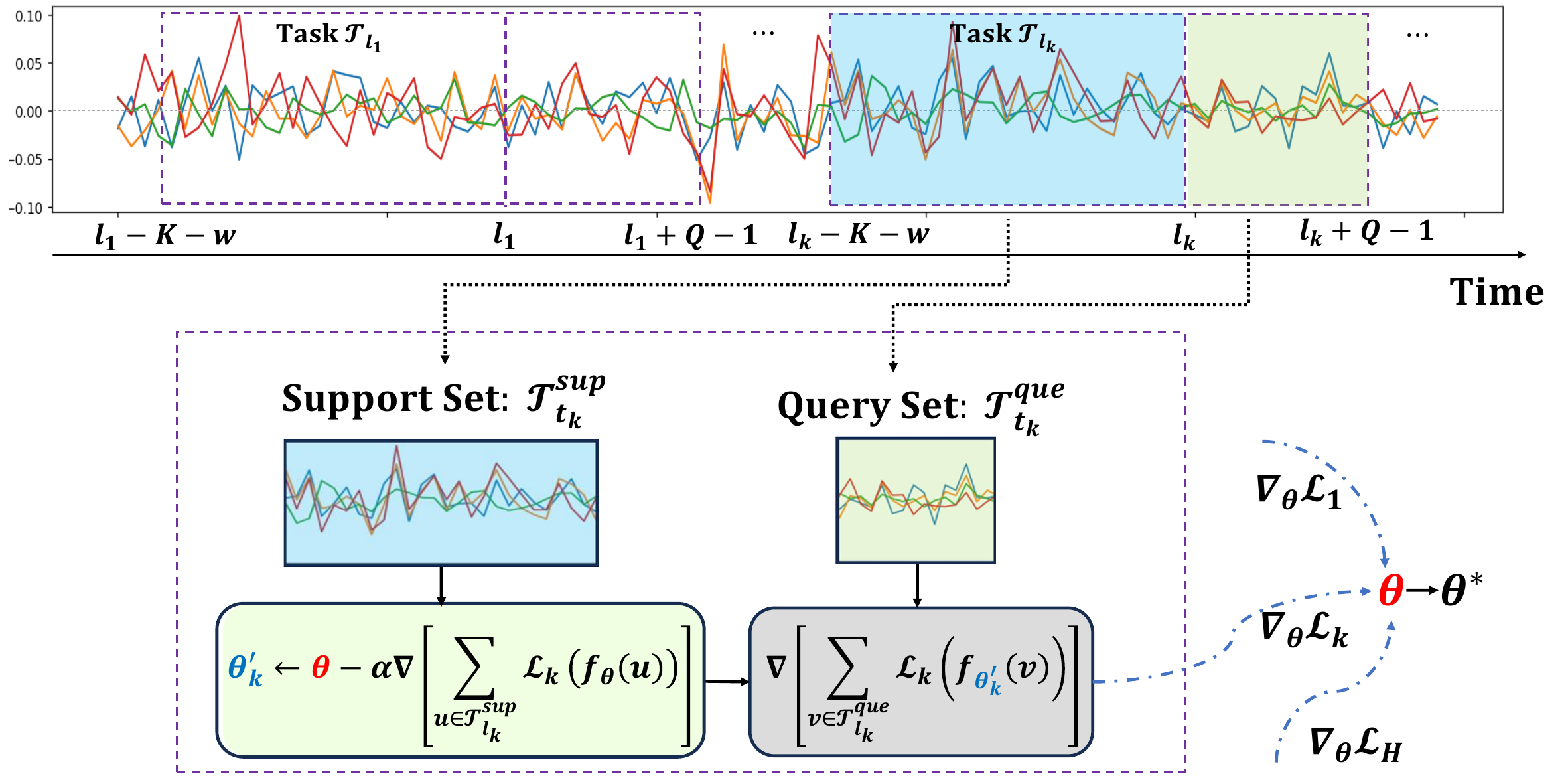}
    \caption{Meta-training framework for online portfolio selection}
    \label{fig:meta—traning framework}
\end{figure}

\textbf{Meta Training.}
Assuming that we have a set of $T^{train}$ training tasks, each consisting of a support set with $K$ samples and a query set with $Q$ samples, which may come from the same market or different markets. 
The training tasks are generated from the original historical return rate data of the $M$ candidate policies by using a rolling window approach. That is, we generate a task in a time window with $K+Q+w$ continue trading days and rolling the time window forward to generate several successive training tasks with overlapping time periods. Please refer to Figure \ref{fig:meta—traning framework} for generation of the training tasks. 

We represent the training tasks as:
\begin{align}
\mathcal{C} = \left\{\mathcal{T}_l = \mathcal{T}_l^{\text{sup}} \cup \mathcal{T}_l^{\text{que}},
\quad l = 1,\dots,
T^{train}\right\},
\end{align}
where
$
\mathcal{T}_l^{\text{sup}} = \left\{\left[\boldsymbol{r}_{t-w}, \dots, \boldsymbol{r}_{t-1}\right]; \hat{\boldsymbol{\omega}}_{t}\right\}_{t=l-K}^{l-1},
$
and
$
\mathcal{T}_l^{\text{que}} = \left\{\left[\boldsymbol{r}_{t-w}, \dots, \boldsymbol{r}_{t-1}\right]; \hat{\boldsymbol{\omega}}_{t}\right\}_{t=l}^{l+Q-1}.
$

We randomly select $H$ tasks $\mathcal{B}=\left\{\mathcal{T}_{l_1},\dots,\mathcal{T}_{l_H}\right\}$ from 
$\mathcal{C}$ for the meta-training. For each selected task $\mathcal{T}_{l_k}$, the support set $\mathcal{T}_{l_k}^{\text{sup}}$ is used to update the 
base parameters $\theta'_k$
for each task, serving as the adaptation process. 
The query set $\mathcal{T}_{l_k}^{\text{que}}$ is used to update the meta-parameter $\theta^*$, helping to find an effective initial parameter that enable the model to quickly adapt to new tasks. The process is described as follows:
\begin{align}
\theta^*=&\operatorname*{argmin}_{\theta}\sum_{k=1}^{H} \mathcal{L}_{k}(f_{\theta_k^{'}}(\mathcal{T}_{l_k}^{que}))\label{eq:OUR}
\\
    &\text{s.t.} \quad\theta_k^{'}=\operatorname*{argmin}_{\theta} \mathcal{L}_k(f_{\theta}(\mathcal{T}_{l_k}^{sup}))
\label{eq:INR}
\end{align}
 where 
 \begin{align}
\mathcal{L}_{k}(f_{\theta_k^{'}}(\mathcal{T}_{l_k}^{que}))&=\sum_{(\boldsymbol{r},\boldsymbol{\hat{\omega}})\in \mathcal{T}_{l_k}^{que}}\mathcal{L}(f_{\theta_k^{'}}(\boldsymbol{r}),\boldsymbol{\hat{\omega}}),\\
\mathcal{L}_k(f_\theta(\mathcal{T}_{l_k}^{sup}))&=\sum_{(\boldsymbol{r},\boldsymbol{\hat{\omega}})\in \mathcal{T}_{l_k}^{sup}}\mathcal{L}(f_\theta(\boldsymbol{r}),\boldsymbol{\hat{\omega}}).
\label{equ:loss_sup}
\end{align}
The optimization of base parameters for each sub-task
\eqref{eq:INR} is performed individually for each task $\mathcal{T}_{l_k}$. It involves updating the base parameters through gradient descent (GD) on the support set $\mathcal{T}_{l_k}^{sup}$:
$$\theta_k^{'}\leftarrow\theta - \alpha \nabla \mathcal{L}_k(f_{\theta}(\mathcal{T}_{l_k}^{sup})),$$ where $\alpha$ is the sub-task learning rate. 
The meta-optimization (\ref{eq:OUR}) updates the meta-parameter $\theta$ by measuring the overall performance of the base parameters $\theta_k'$ based on $\theta$ on the query set $\mathcal{T}_{l_k}^{que}$. 
The meta-parameters is update by the stochastic gradient descent (SGD) on the $H$ randomly selected  tasks:
$$\theta^* \leftarrow \theta - \beta \nabla \left[\sum_{k=1}^{H} \mathcal{L}_k(f_{{\theta_k^{'}}}(\mathcal{T}_{l_k}^{que}))\right],$$ where $\beta$ is the meta-learning rate. We illustrate this process in Figure \ref{fig:meta—traning framework}.

In the proposed meta-learning framework, we decompose a long-term portfolio selection problem into multiple short-term tasks, which allows us to treat the global non-stationary time series as locally stationary. As a result, meta-learning helps mitigate the bias caused by changes in temporal distributions and enhances the model's generalizability.

\textbf{Meta Testing.}
After completing the meta-training process, we obtain the initial model parameters \(\theta^*\). We then use the support set of the test task \(\{\mathcal{T}_l^{sup}\}\) to perform gradient descent (GD) on \(\theta^*\), rapidly acquiring model parameters tailored to this specific task, denoted as \(\theta^{new}_l\). Next, by feeding the features into the model \(f_{\theta^{new}_l}\), we obtain the weight allocation for the \(M\) candidate policies, which allows us to determine the final investment vector.

refWe refer to this investment portfolio management algorithm, which is based on meta-learning and mixture policies learning, as the Meta-LMPS-Online algorithm. The algorithm flow is presented in Algorithm \ref{Meta-LMPS-Online}.


\begin{algorithm}[ht]
    \caption{Meta-LMPS-Online
    \label{Meta-LMPS-Online}}
        \renewcommand{\algorithmicrequire}{\textbf{Input:}}
        \renewcommand{\algorithmicensure}{\textbf{Output:}}
    
    \begin{algorithmic}[1]
    \REQUIRE candidate policies $\{\pi^j (\boldsymbol{b}|\boldsymbol{x})\}_{j=1}^M$, historical price vector sequence $\boldsymbol{x}$
        \ENSURE Model parameters $\theta_l^{new}$ suitable for test task $\mathcal{T}_l$
        \STATE Construct the set of training tasks set $\mathcal{C}=\{\mathcal{T}\} $ {using candidate policies and historical price vector sequence}
        \FOR{each epoch}
        \STATE Randomly select $H$ tasks $\mathcal{B}:\{\mathcal{T}_{l_1},\dots,\mathcal{T}_{l_H}\}$ from $\mathcal{C}$
        \FOR{$\mathcal{T}_{l_k}$}
         \STATE Evaluate $\nabla_\theta \mathcal{L}_{k}(f_\theta(\mathcal{T}_{l_k}^{sup}))$ using $\mathcal{T}_{l_k}^{sup}$ and $\mathcal{L}_k$ in equation \eqref{equ:loss_sup}  
        \STATE Compute base parameters with gradient descent: $\theta_k'\leftarrow\theta-\alpha \nabla_\theta \mathcal{L}_k(f_\theta(\mathcal{T}_{l_k}^{sup}))$
        \ENDFOR
        \STATE Update 
        $\theta\leftarrow \theta-\beta\nabla_\theta\sum_{k=1}^H \mathcal{L}_k(f_{\theta_k^{'}}(\mathcal{T}_{l_k}^{que}))$
        \ENDFOR
        \STATE Acquire optimized meta-parameter $\theta^*\leftarrow\theta$
        \FOR{trading period $l$}
        \STATE Set the initial parameter $\theta_0 = \theta^*$
        \STATE Construct the support set $\mathcal{T}_l^{sup}$ using price vectors
        \STATE Perform gradient descent on $\mathcal{T}_{l}^{sup}$ to quickly obtain $\theta_l^{new} = \theta_0 - \alpha \nabla_\theta \mathcal{L}_l(f_{\theta_0}(\mathcal{T}_l^{sup}))$
        \ENDFOR
    \end{algorithmic}
\end{algorithm}

\begin{figure}[htp]
    \centering
        \includegraphics[width=1.1\textwidth]{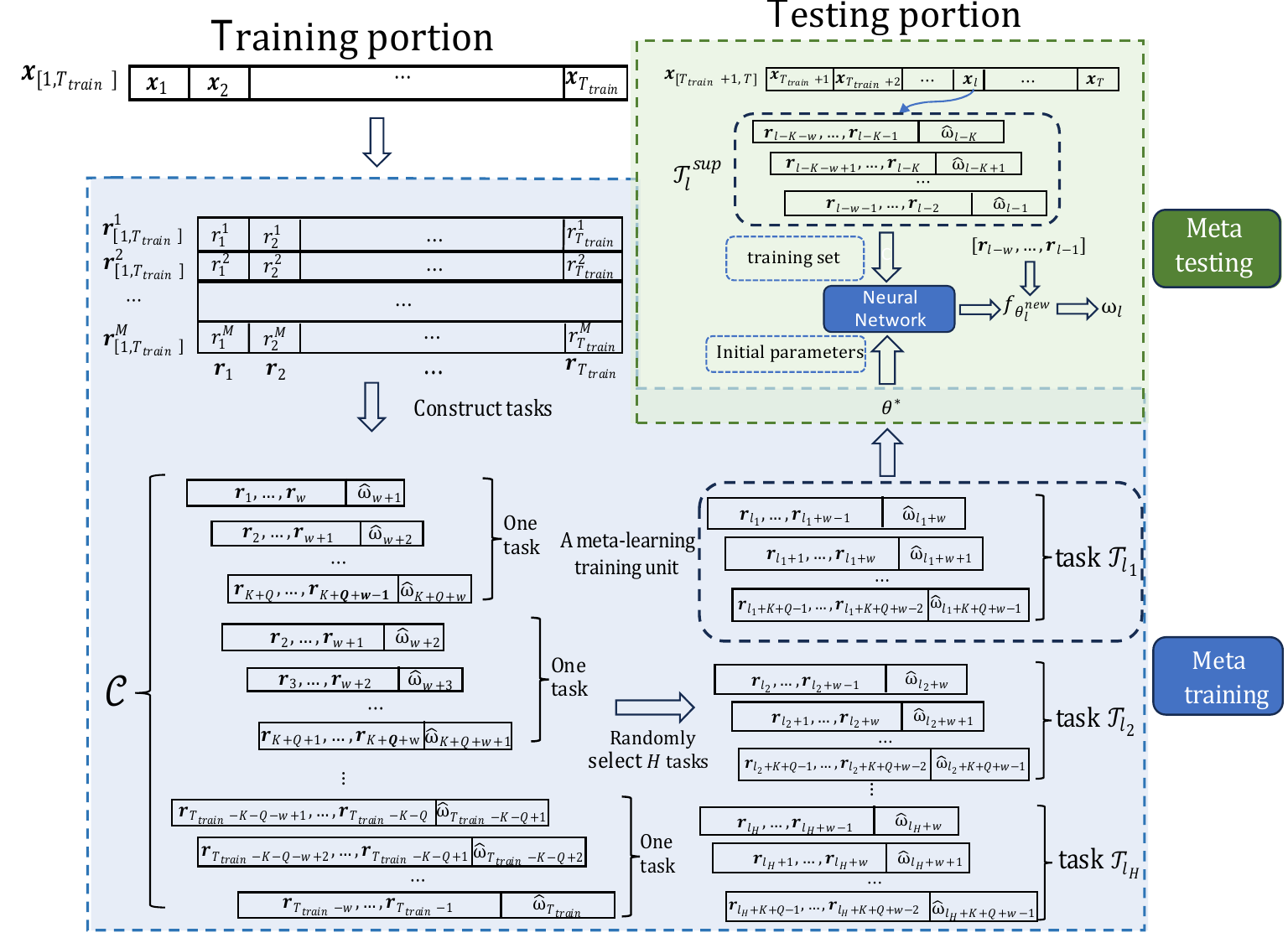}
        \caption{The workflow of LMPS-SMO algorithm}
        \label{fig:The_workflow_of_LMPS_SMO_algorithm}
  \end{figure}

\begin{figure}[htp]
        \centering
        \includegraphics[width=1.1\textwidth]{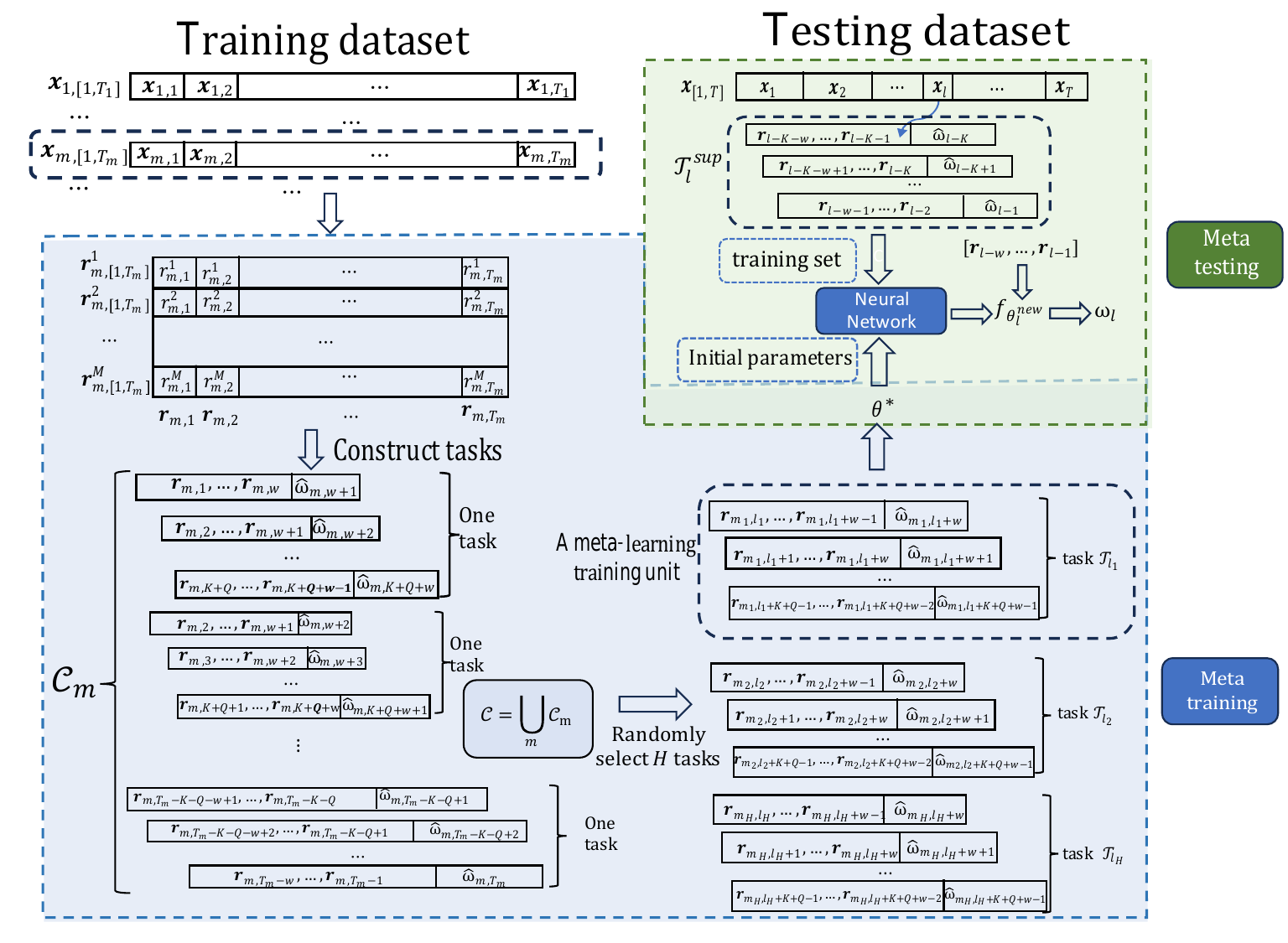}
        \caption{The workflow of LMPS-CMO algorithm}
        \label{fig:The_workflow_of_LMPS_CMO_algorithm}
\end{figure}

\textbf{Single-market and Cross-market algorithm.}
The Meta-LMPS-Online model features two sub-case algorithms for two types of task data: 
the Single-market Adapted Meta-LMPS-Online Algorithm (LMPS-SMO) and the Cross-Market Learning Meta-LMPS-Online Algorithm (LMPS-CMO). LMPS-SMO uses training data and test data from the same market, allowing for swithing of market regimes.
While LMPS-CMO utilizes training data and test data from multiple different markets, demonstrating broader generalization. Unlike other methods that are limited 
to a single market due to the data-consistency issue, our approach facilitates the transfer of learning experiences across different markets, offering greater flexibility and adaptability. The detailed algorithms and illustrations for LMPS-SMO and LMPS-CMO are provided in Algorithm \ref{Algorithm Overview in a single market} and Figure \ref{fig:The_workflow_of_LMPS_SMO_algorithm} for the former, and in Algorithm \ref{Algorithm Overview under mutiple datasets} and Figure \ref{fig:The_workflow_of_LMPS_CMO_algorithm} for the latter.
\begin{algorithm}[p]
\begin{footnotesize} 
\caption{\footnotesize{Single-Market Adapted Meta-LMPS-Online Algorithm(LMPS-SMO)}}
\label{Algorithm Overview in a single market}
\renewcommand{\algorithmicrequire}{\textbf{Input:}}
\renewcommand{\algorithmicensure}{\textbf{Output:}}

\begin{algorithmic}[1]
\REQUIRE ~~ $\boldsymbol{x}_{[1,T_{train}]}$: Training portion of price relative vector sequences
    \REQUIRE ~~ $\boldsymbol{x}_{[T_{train}+1,T]}$: Testing portion of price relative vector sequences
    \REQUIRE ~~ $M$: Number of candidate policies
    \REQUIRE ~~ $K,Q$: Number samples in support set, Number samples in query set
    \REQUIRE ~~ $H$: Number of batch training tasks
    \ENSURE ~~ $S_{[T_{train}+1,T]}$: Cumulative wealth sequences in the testing portion\\
    \vspace{0.5\baselineskip}
    //Data Preparation and Data Preprocessing
    \STATE Based on \(\boldsymbol{x}_{[1,T_{\text{train}}]}\), use a clustering method to identify $M$ historically high-performing policies with low similarity as candidate strategies. 
    \STATE Obtain the return sequences of $M$ candidate policies $\boldsymbol{r}_{[1,T_{train}]}^j,j=1,\dots,M$ on training portion $\boldsymbol{x}_{[1,T_{\text{train}}]}$, and construct training task set $\mathcal{C}=\{\mathcal{T}\}$

    \vspace{0.5\baselineskip}
    //Meta Training
    \STATE Randomly initialize $\theta$\\
    \FOR{each epoch}
    \STATE Randomly select $H$ tasks $\mathcal{B}=\{\mathcal{T}_{l_k}\}_{k=1}^H$ from $\mathcal{C}$ \\
    \FOR{$\mathcal{T}_{l_k}$}
     \STATE Evaluate $\nabla_\theta \mathcal{L}_k(f_\theta(\mathcal{T}_{l_k}^{sup}))$ with respect to $K$ examples
    \STATE Compute adapted parameters with gradient descent: $\theta_k'=\theta-\alpha \nabla_\theta \mathcal{L}_k(f_\theta(\mathcal{T}_{l_k}^{sup}))$
    \ENDFOR
    \STATE Update $\theta\leftarrow \theta-\beta\nabla_\theta\sum_{k=1}^H \mathcal{L}_{\mathcal{T}_{k}}(f_{\theta_k^{'}}(\mathcal{T}_{l_k}^{que}))$
    \ENDFOR\\
    \STATE Update $\theta^*\leftarrow\theta$

    \vspace{0.5\baselineskip}
    //Meta Testing 
    \STATE Initialization: $S_{T_{train}+1}=1$
    \FOR{$l=T_{train}+1,\dots,T$}
    \STATE Set the initial parameter $\theta_0=\theta^*$
    \STATE Calculate the returns of $M$ candidate policies in period $(t-1)$ to obtain the support set for the current meta-test task $\mathcal{T}_{l}^{sup}=\{[\boldsymbol{r}_{l+t-w},\dots,\boldsymbol{r}_{l+t-1}];\hat{\boldsymbol{\omega}}_{l+t}\}_{i=-K}^{-1}$
    \STATE Perform gradient descent on $\mathcal{T}_{l}^{sup}$ to quickly obtain $\theta^{new} = \theta_0-\alpha \nabla_\theta \mathcal{L}_t(f_{\theta_0}(\mathcal{T}_t^{sup}))$
    \STATE Obtain the mixture weights in period $t$: $\boldsymbol{\omega}_t=f_{\theta^{new}_t}(\boldsymbol{r}_{t-w},\dots,\boldsymbol{r}_{t-1})$
    \STATE Acquire the portfolios of $M$ candidate policies in period $l$ ($\boldsymbol{b}_{l}^1,\dots,\boldsymbol{b}_{l}^M)$, and set the investment portfolio as $\boldsymbol{b}_{l}=\sum_{j=1}^M \omega_{l}^j \boldsymbol{b}_{l}^j$
    \STATE The market reveals actual price relative vectors $\boldsymbol{x}_l$
    \STATE Update the cumulative wealth as $S_l=S_{l-1}\cdot \boldsymbol{b}_l^\top \boldsymbol{x}_l(1-\delta\Vert \boldsymbol{b}_l-\widetilde{\boldsymbol{b}}_{l-1}\Vert_1)$
    \ENDFOR
\end{algorithmic}
\end{footnotesize}
\end{algorithm}

\begin{algorithm}[p]
\begin{footnotesize} 
\caption{\footnotesize{Cross-Market Learning LMPS-Online Algorithm(LMPS-CMO)}
    \label{Algorithm Overview under mutiple datasets}}
    \renewcommand{\algorithmicrequire}{\textbf{Input:}}
    \renewcommand{\algorithmicensure}{\textbf{Output:}}

\begin{algorithmic}[1]
    \REQUIRE ~~ $\boldsymbol{x_{1,[1,T_1]}},\dots,\boldsymbol{x_{m,[1,T_m]}},\dots$: Relative price sequences of training datasets
        \REQUIRE ~~ $\boldsymbol{x}_{[1,T]}$: Relative price sequences of testing dataset
        \REQUIRE ~~ $M$: Number of candidate policies
        \REQUIRE ~~ $K,Q$: Number samples in support set, Number samples in query set
        \REQUIRE ~~ $H$: Number of batch training tasks
        \ENSURE ~~ $S_{[K+w,T]}$: Cumulative wealth sequences in the testing portion\\
        
        //Data Preparation and Data Preprocessing
        \STATE Based on training datasets, use a clustering method to identify $M$ historically high-performing policies with low similarity as candidate strategies. 
        \STATE Initialize $\mathcal{C}=\{\}$
        \FOR{each $m$}
        \STATE Obtain the return sequences of $M$ candidate policies $\boldsymbol{r}_{m,[1,T_m]}^j,j=1,\dots,M$ on training portion $\boldsymbol{x}_{m,[1,T_m]}$, and construct task set $\mathcal{C}_m=\{\mathcal{T}\}$
        \STATE Update $\mathcal{C}=\mathcal{C}\bigcup\mathcal{C}_m$
        \ENDFOR
        
        \vspace{0.5\baselineskip}
        //Meta Training
        \STATE Randomly initialize $\theta$\\
        \FOR{each epoch}
        \STATE Randomly select $H$ tasks $\mathcal{B}=\{\mathcal{T}_{m_k,l_k}\}_{k=1}^H$ from $\mathcal{C}$ \\
        \FOR{$\mathcal{T}_{m_k,l_k}$}
         \STATE Evaluate $\nabla_\theta \mathcal{L}_k(f_\theta(\mathcal{T}_{m_k,l_k}^{sup}))$ with respect to $K$ examples
        \STATE Compute adapted parameters with gradient descent: $\theta_k'=\theta-\alpha \nabla_\theta \mathcal{L}_k(f_\theta(\mathcal{T}_{m_k,l_k}^{sup}))$
        \ENDFOR
        \STATE Update $\theta\leftarrow \theta-\beta\nabla_\theta\sum_{k=1}^H \mathcal{L}_{\mathcal{T}_{k}}(f_{\theta_k^{'}}(\mathcal{T}_{m_k,l_k}^{que}))$
        \ENDFOR\\
        \STATE Update $\theta^*\leftarrow\theta$

        \vspace{0.5\baselineskip}
        //Meta Testing
        \STATE Initialization:$S_{K+w}=1$
        \FOR{$l=K+w+1,\dots,T$}
        \STATE Set the initial parameter $\theta_0=\theta^*$
        \STATE Calculate the returns of $M$ candidate policies in period $(t-1)$ to obtain the support set for the current meta-test task $\mathcal{T}_{l}^{sup}=\{[\boldsymbol{r}_{l+t-w},\dots,\boldsymbol{r}_{l+t-1}];\hat{\boldsymbol{\omega}}_{l+t}\}_{i=-K}^{-1}$
        \STATE Perform gradient descent on $\mathcal{T}_{l}^{sup}$ to quickly obtain $\theta^{new} = \theta_0-\alpha \nabla_\theta \mathcal{L}_t(f_{\theta_0}(\mathcal{T}_t^{sup}))$
        \STATE Obtain the mixture weights in period $t$: $\boldsymbol{\omega}_t=f_{\theta^{new}_t}(\boldsymbol{r}_{t-w},\dots,\boldsymbol{r}_{t-1})$
        \STATE Acquire the portfolios of $M$ candidate policies in period $l$ ($\boldsymbol{b}_{l}^1,\dots,\boldsymbol{b}_{l}^M)$, and set the investment portfolio as $\boldsymbol{b}_{l}=\sum_{j=1}^M \omega_{l}^j \boldsymbol{b}_{l}^j$
        \STATE The market reveals actual price relative vectors $\boldsymbol{x}_l$
        \STATE Update the cumulative wealth as $S_l=S_{l-1}\cdot \boldsymbol{b}_l^\top \boldsymbol{x}_l(1-\delta\Vert \boldsymbol{b}_l-\widetilde{\boldsymbol{b}}_{l-1}\Vert_1)$
        \ENDFOR
\end{algorithmic}
\end{footnotesize}
\end{algorithm}

\newpage

\section{Experiments}
\subsection{Experimental Setting}
\textbf{Datasets.}
Experiments are conducted on six real-world datasets of daily stock returns, namely NYSE(O), NYSE(N), SP500, DJIA, TSE, and MSCI. The detailed information for each dataset is shown in Table \ref{tab:Datasets Information}. The datasets were downloaded from the \href{https://github.com/OLPS/OLPS}{OLPS GitHub repository}.

\begin{table}[H]
	\centering
	\caption{Datasets Information}
	\label{tab:Datasets Information}{%
		\begin{tabular}{ccccccc}
			\hline
			Dataset Name & Region & Time Interval & Number of Assets & Frequency  \\ \hline
			\textbf{NYSE(O)} & America & 5651 & 36 & Daily \\
			\textbf{NYSE(N)} & America & 6431 & 23 & Daily \\
			\textbf{SP500} & America & 1276 & 25 & Daily \\
                \textbf{DJIA}  & America & 507 & 30 & Daily \\
   			\textbf{TSE}  & Canada & 1259 & 88 & Daily \\
			\textbf{MSCI} & China & 1043 & 24 & Daily \\
			\hline		
		\end{tabular}%
	}
\end{table}

\textbf{Baseline.}
We compare our method with 14 traditional policies listed in Table \ref{tab:policy}. 

\textbf{Experimental Setup.}
The transaction cost rate for buying and selling assets is set to 0.0005, as referenced from the configuration in \cite{GUO2023}. The number of candidate policies \(M\) is set as 4. 
The historical window size $w$ 
is 5. 
The sample size \(K\) of support set of each task is 10, while the sample size \(Q\) of the query set is 1. 
The epoch of each task is 
20. 
The base-parameter learning rate \(\alpha\) is 0.001, and the meta-learning rate \(\beta\) is 0.0005.

\subsection{Comparison of LMPS-SMO to Baselines}

\begin{figure}[h!]
    \centering
    \includegraphics[width=\textwidth]{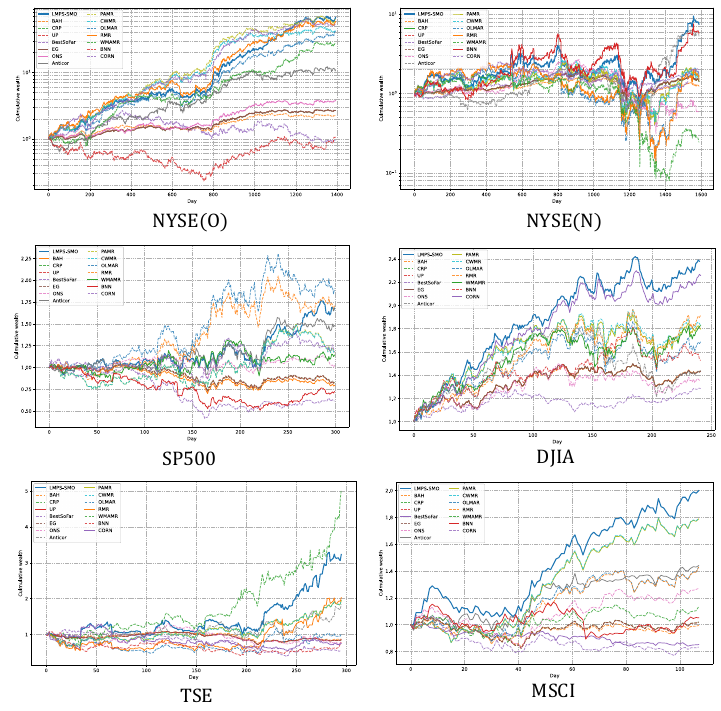}
    \caption{Single-market cumulative return curves of different strategies in $6$ datasets. }
    \label{fig:Cumulative wealth for strategies under $6$ datasets}
\end{figure}

\textbf{Single-market test.}
For each dataset, We utilize the first three quarters of the data for meta-training purposes, while the remaining one quarter of the data is employed for meta-testing.
Table \ref{table:Comparison of Returns and Risks on 6 datasets} and Figure \ref{fig:Cumulative wealth for strategies under $6$ datasets} summarize the comparison results of LMPS-SMO and the baseline method.

\begin{table*}[htp]
\centering
\caption{Single-market performance of different strategies in 6 datasets}
\label{table:Comparison of Returns and Risks on 6 datasets}
\footnotesize
\begin{tabular}{ccccccc}
\hline
 Datasets& \multicolumn{2}{c}{NYSE(O)} & \multicolumn{2}{c}{NYSE(N)} & \multicolumn{2}{c}{SP500}  
 \\ \hline
Strategy & \begin{tabular}[c]{@{}c@{}}Cumulative\\ Wealth \end{tabular} & \begin{tabular}[c]{@{}c@{}}Sharpe \\Ratio\end{tabular} & \begin{tabular}[c]{@{}c@{}}Cumulative\\ Wealth\end{tabular} & \begin{tabular}[c]{@{}c@{}}Sharpe \\Ratio\end{tabular} & \begin{tabular}[c]{@{}c@{}}Cumulative\\ Wealth\end{tabular} & \begin{tabular}[c]{@{}c@{}}Sharpe \\Ratio\end{tabular}  
\\ \hline
\textbf{LMPS-SMO} & {\color{red}\textbf{67.44}} & 19.05 & {\color{red}\textbf{7.39}} & {\color{red}\textbf{3.55}} & 1.68 & 17.41  \\
BAH & 2.24 & 6.30 & 1.25 & 1.08 & 0.78$^c$ & -10.44$^c$  \\
CRP & 2.68 & 8.93 & 1.43$^c$ & 1.46$^c$ & 0.82 & -8.55 \\
UP & 2.68 & 8.93 & 1.44 & 1.50 & 0.82 & -8.59  \\
BestSoFar & 0.94 & -0.20 & 1.77 & 2.61 & 0.59 & -11.41 \\
EG & 2.68$^c$ & 8.90$^c$ & 1.44 & 1.51 & 0.82 & -8.62  \\
ONS & 3.84$^c$ & 9.09$^c$ & 0.59 & -0.74 & 1.02 & 0.78 \\
Anticor & 11.00$^c$ & 9.46$^c$ & 5.32 & 3.16 & 1.48$^c$ & 12.14$^c$  \\
PAMR & 62.33 & {\color{red}\textbf{19.51}} & 1.59$^c$ & 0.72$^c$ & 1.15 & 3.79  \\
CWMR & 42.45 & 17.46 & 1.49 & 0.62 & 1.14 & 3.59  \\
OLMAR & 41.34 & 14.60 & 1.61 & 0.58 & {\color{red}\textbf{1.86}} & {\color{red}\textbf{20.26}} \\
RMR & 56.53$^c$ & 16.38$^c$ & 1.46$^c$ & 0.47$^c$ & 1.67 & 15.97 \\
WMAMR & 28.31 & 13.28 & 0.24 & -1.55 & 1.13$^c$ & 3.21$^c$ \\
BNN & 1.01 & 0.02 & 5.95$^c$ & 2.73$^c$ & 0.72$^c$ & -8.21$^c$ \\
CORN & 53.00 & 18.24 & 1.83 & 1.43 & 1.13 & 3.59  \\
\hline
 Datasets& 
 \multicolumn{2}{c}{DJIA} & \multicolumn{2}{c}{TSE} & \multicolumn{2}{c}{MSCI} \\ \hline
Strategy & \begin{tabular}[c]{@{}c@{}}Cumulative\\ Wealth \end{tabular} & \begin{tabular}[c]{@{}c@{}}Sharpe \\Ratio\end{tabular} & \begin{tabular}[c]{@{}c@{}}Cumulative\\ Wealth\end{tabular} & \begin{tabular}[c]{@{}c@{}}Sharpe \\Ratio\end{tabular} & \begin{tabular}[c]{@{}c@{}}Cumulative\\ Wealth\end{tabular} & \begin{tabular}[c]{@{}c@{}}Sharpe \\Ratio\end{tabular} \\
\hline
\textbf{LMPS-SMO} & {\color{red}\textbf{2.00}} & {\color{red}\textbf{208.71}} & 3.23 & 34.42 & {\color{red}\textbf{2.38}} & {\color{red}\textbf{90.77}}\\
BAH & 1.02 & 3.86 & 0.89 & -7.79 & 1.37 & 31.26\\
CRP & 1.05 & 8.88 & 0.86 & -9.23 & 1.36 & 29.88\\
UP & 1.05 & 8.66 & 0.86$^c$ & -9.23$^c$ & 1.36 & 30.02\\
BestSoFar  & 0.86$^c$ & -34.19$^c$ & 0.56 & -10.35 & 1.27 & 25.74\\
EG & 1.05 & 8.62 & 0.86 & -9.31 & 1.36$^c$ & 29.99$^c$\\
ONS  & 1.30 & 49.49 & 0.81 & -4.68 & 1.25 & 15.03\\
Anticor & 1.50$^c$ & 80.3$^c$ & 1.86 & 12.29 & 1.42 & 20.64\\
PAMR & 1.80$^c$ & 171.64$^c$ & 1.93$^c$ & 17.95$^c$ & 1.69$^c$ & 35.04$^c$\\
CWMR& 1.81 & 172.48 & 1.94 & 18.15 & 1.72 & 36.26\\
OLMAR & 1.45 & 73.95 & 1.00 & 0.05 & 1.59 & 28.40\\
RMR & 1.44 & 71.74 & 2.06$^c$ & 12.81$^c$ & 1.80 & 39.18 \\
WMAMR & 1.15 & 21.33 & {\color{red}\textbf{4.93}} & {\color{red}\textbf{55.29}} & 1.75$^c$ & 37.45$^c$ \\
BNN & 1.08$^c$ & 11.05$^c$ & 0.80 & -3.8 & 1.39 & 20.92 \\
CORN& 0.85 & -20.21 & 0.76$^c$ & -4.76$^c$ & 2.17$^c$ & 87.63$^c$\\
\hline
\multicolumn{6}{l}{$\cdot^c$: Candidate policies selected for each dataset.}
\end{tabular}
\end{table*}

From Table \ref{table:Comparison of Returns and Risks on 6 datasets}, it can be observed that the LMPS-SMO method outperforms the candidate policies in terms of cumulative wealth and Sharpe ratio across all datasets. As shown in Figure \ref{fig:Cumulative wealth for strategies under $6$ datasets}, the LMPS-SMO method generally tends to learn one or two candidate policies while also being able to adjust policy allocations in a timely manner, resulting in its overall performance surpassing that of the individual learned candidate policies. 
Although some traditional policies perform well across these datasets, the LMPS-SMO method consistently ranks high in terms of return and risk control metrics. 
By learning from the FTL and pattern-matching strategies, the LMPS-SMO strategy achieves higher returns, while learning from the benchmark and FTW strategies steers it towards stability, thereby effectively reducing risk.

\begin{figure*}[ht]
    \centering
    \includegraphics[width=\textwidth]{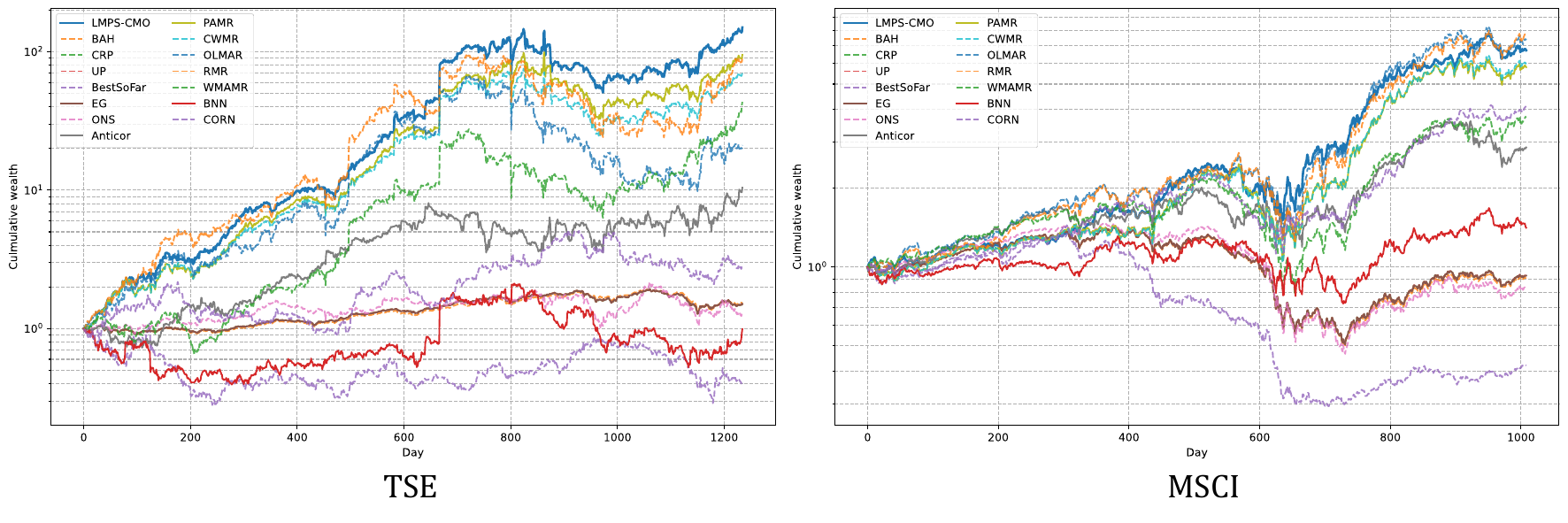}
    \caption{Cross-market cumulative return curves of different strategies in TSE and MSCI datasets}
    \label{fig:Cumulative wealth for strategies under $2$ datasets}
\end{figure*}

\textbf{Cross-Market test:} We use four datasets (NYSE(O), NYSE(N), DJIA and SP500) to construct training tasks and select candidate policies. We then evaluate LMPS-CMO method using the TSE and MSCI datasets as the test set for online portfolio management. The results are presented in Table \ref{tab:cross_market} and Figure \ref{fig:Cumulative wealth for strategies under $2$ datasets}. 

\begin{table}[htb]
    \caption{Cross-market performance of different policies in TSE and MSCI datasets}
    \centering
    \scalebox{0.85}[0.85]{
    \label{tab:cross_market}
    \begin{tabular}{ccccc}
\hline
Datasets & \multicolumn{2}{c}{TSE} & \multicolumn{2}{c}{MSCI} \\
\hline
Strategy & \begin{tabular}[c]{@{}c@{}}Culmulative \\ Wealth\end{tabular} & \begin{tabular}[c]{@{}c@{}}Sharpe\\ Ratio\end{tabular} & \begin{tabular}[c]{@{}c@{}}Culmulative \\ Wealth\end{tabular} & \begin{tabular}[c]{@{}c@{}}Sharpe\\ Ratio\end{tabular} \\
\hline
\textbf{LMPS-CMO} & {\color{red} \textbf{149.57}} & {\color{red} \textbf{17.54}} & 6.73 & 13.30 \\
BAH & 1.55 & 5.10 & 0.93 & -0.55 \\
CRP & 1.51 & 4.76 & 0.95 & -0.40 \\
UP & 1.51 & 4.78 & 0.95 & -0.39 \\
BestSoFar & 0.42 & -2.04 & 0.43 & -6.07 \\
EG & 1.51$^c$ & 4.76$^c$ & 0.95$^c$ & -0.40$^c$ \\
ONS & 1.34 & 1.28 & 0.85 & -1.00 \\
Anticor & 10.40$^c$ & 6.81$^c$ & 2.93$^c$ & 6.51$^c$ \\
PAMR & 95.07$^c$ & 15.01$^c$ & 5.91$^c$ & 11.86$^c$ \\
CWMR & 71.49 & 13.76 & 6.07 & 12.04 \\
OLMAR & 21.36 & 6.50 & 7.76 & 13.10 \\
RMR & 90.24 & 11.38 & {\color{red} \textbf{8.09}} & {\color{red} \textbf{13.51}} \\
WMAMR & 43.32 & 9.55 & 3.92 & 8.03 \\
BNN & 1.01$^c$ & 0.02$^c$ & 1.44$^c$ & 2.21$^c$ \\
CORN & 2.92 & 3.03 & 4.22 & 10.74\\
\hline
\multicolumn{5}{l}{$\cdot^c$: Candidate policies selected for each dataset.}
\end{tabular}
}
\end{table}

The results demonstrate that within a framework of cross-learning from multiple datasets, the Meta-LMPS-online algorithm maintains strong performance. Whether on the TSE or MSCI datasets, the LMPS-CMO strategy showcases higher cumulative wealth and Sharpe ratios compared to the four candidate policies it learned. Moreover, it manages to contain volatility and maximum drawdown within a narrower range. This underscores the algorithm's robust learning capacity and exceptional adaptability to various environments.

\subsection{Ablation studies}
\begin{figure*}[h!]
    \centering
    \includegraphics[width=\textwidth]{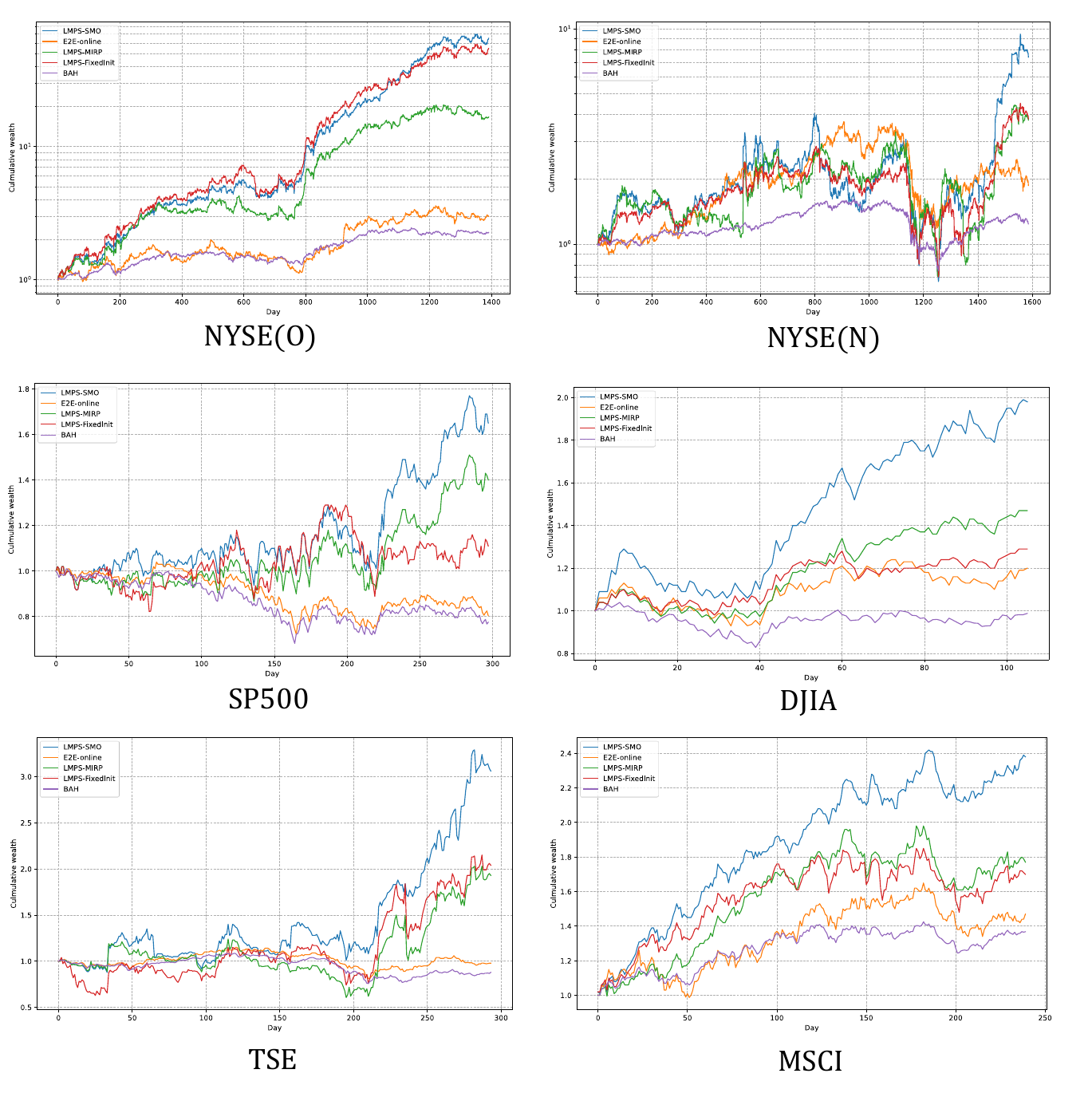}
    \caption{Ablation study of the LMPS-SMO method in $6$ datasets}
    \label{fig:Ablation studies1}
\end{figure*}

We carry out a series of ablation studies
to examine the affect of meta-training and meta-testing  on the performance of the Meta-LMPS-online algorithm. For this purpose, we designed the following three controlled experiences
to compare against the completed
Meta-LMPS-online algorithm. 
\begin{itemize}
    \item E2E-online: Without using the mixture policies learning framework, the network takes the return sequences of assets as input and the allocation ratios for each asset as output. Apart from this, other settings remain the same as the LMPS-online algorithm.
    \item LMPS-MIRP: Without using meta-training to generate initial parameters for meta-testing, it randomly generates initial parameters and updates these parameters during the meta-testing phase.
    \item LMPS-FixedInit: During testing phase, it consistently uses the initial parameters generated in the meta-training phase to determine policy allocations, without any updating in the meta-testing period.
\end{itemize}
\textbf{Single-market test:}
The final results are displayed in Table \ref{tab:single_stage_object2} and Figure \ref{fig:Ablation studies1}. We include the BAH strategy for comparison to observe the overall market state.

The results show that the LMPS-SMO strategy outperforms others in cumulative returns across 6 datasets. The E2E-online strategy follows market trends but generally performs worse than the mixture policies learning algorithm, demonstrating the benefits of the mixture policies approach. LMPS-FixedInit tends to perform well for specific datasets but lacks adaptability across different market conditions, highlighting the importance of meta-testing. LMPS-MIRP also performs worse than LMPS-SMO, emphasizing the need for meta-training.

\begin{table*}[htpb]
\centering
\caption{Ablation study of the LMPS-SMO method in 6 datasets}
\label{tab:single_stage_object2}
\resizebox{\textwidth}{!}
{
\renewcommand{\arraystretch}{1.1} 
\begin{tabular}{cccccccccccccc}
\hline
Datasets & Strategy & \begin{tabular}[c]{@{}c@{}}Culmulative\\ Wealth\end{tabular} & \begin{tabular}[c]{@{}c@{}}Annualized\\ Return\end{tabular} & \begin{tabular}[c]{@{}c@{}}Annualized\\ Volatility\end{tabular} & \begin{tabular}[c]{@{}c@{}}Maximum\\ Drawdown\end{tabular} & \begin{tabular}[c]{@{}c@{}}Sharpe\\ Ratio\end{tabular} & 
\\
\hline
 & \textbf{LMPS-SMO} & {\color{red} \textbf{67.44}} & {\color{red} \textbf{111.76}\%} & 5.93\% & 27.43\% & {\color{red} \textbf{19.05}} &   \\
 & E2E-online & 2.99 & 21.74\% & 4.64\% & 44.00\% & 4.68 &  \\
 & LMPS-MIRP & 16.70 & 65.80\% & 5.73\% & 35.05\% & 11.48 &  \\
 & LMPS-FixedInit & 54.70 & 105.18\% & 6.50\% & 43.60\% & 16.17 & \\
\multirow{-5}{*}{NYSE(O)} & BAH & 2.26 & 15.80\% & {\color{red} \textbf{2.47}\%} & {\color{red} \textbf{20.37}\%} & 6.39 & \\
\hline
 & \textbf{LMPS-SMO} & {\color{red} \textbf{7.39}} & {\color{red} \textbf{37.04}\%} & 10.43\% & 83.48\% & {\color{red} \textbf{3.55}} &   \\
 & E2E-online & 1.87 & 10.36\% & 6.67\% & 68.19\% & 1.55 &  \\
 & LMPS-MIRP & 3.77 & 23.25\% & 12.00\% & 78.67\% & 1.94 &  \\
 & LMPS-FixedInit & 3.82 & 23.51\% & 7.99\% & 75.33\% & 2.94 &   \\
\multirow{-5}{*}{NYSE(N)} & BAH & 1.24 & 3.45\% & {\color{red} \textbf{3.21}\%} & {\color{red} \textbf{53.59}\%} & 1.07 &  \\
\hline
 & \textbf{LMPS-SMO} & {\color{red} \textbf{1.68}} & {\color{red} \textbf{52.21}\%} & 3.12\% & 22.48\% & {\color{red} \textbf{17.41}} &   \\
 & E2E-online & 0.80 & -16.72\% & {\color{red} \textbf{1.79}\%} & 30.48\% & -9.34 &   \\
 & LMPS-MIRP & 1.40 & 32.61\% & 3.04\% & {\color{red} \textbf{22.46}\%} & 10.71 &   \\
 & LMPS-FixedInit & 1.11 & 9.15\% & 3.25\% & 31.16\% & 2.82 &   \\
\multirow{-5}{*}{SP500} & BAH & 0.78 & -19.24\% & 1.84\% & 32.46\% & -10.48 & \\
\hline
& \textbf{LMPS-SMO} & {\color{red} \textbf{2.00}} & {\color{red} \textbf{400.82}\%} & 1.90\% & 17.83\% & {\color{red} \textbf{208.71}} \\
& E2E-online & 1.20 & 53.73\% & 1.61\% & 17.79\% & 33.47  \\
& LMPS-MIRP & 1.47 & 148.09\% & 1.37\% & 14.36\% & 108.18 \\
 & LMPS-FixedInit & 1.29 & 82.32\% & 1.30\% & {\color{red} \textbf{10.91}\%} & 63.21 \\
\multirow{-5}{*}{DJIA} & BAH & 0.99 & -2.87\% & {\color{red} \textbf{1.16}\%} & 20.31\% & -2.47 \\
& \textbf{LMPS-SMO} & {\color{red} \textbf{3.23}} & {\color{red} \textbf{158.84}\%} & 4.91\% & 32.22\% & {\color{red} \textbf{34.42}}\\
& E2E-online & 0.98 & -1.87\% & 1.29\% & {\color{red} \textbf{29.48}\%} & -1.46 \\
& LMPS-MIRP & 1.93 & 74.91\% & 5.17\% & 51.13\% & 14.50\\
& LMPS-FixedInit & 2.04 & 83.35\% & 5.55\% & 39.52\% & 15.01\\
\multirow{-5}{*}{TSE} & BAH & 0.88 & -10.46\% & {\color{red} \textbf{1.25}\%} & 30.02\% & -8.39\\
\hline
& \textbf{LMPS-SMO} & {\color{red} \textbf{2.38}} & {\color{red} \textbf{146.76}\%} & 1.62\% & {\color{red} \textbf{12.40}\%} & {\color{red} \textbf{90.77}}\\
& E2E-online & 1.47 & 49.38\% & 2.29\% & 21.67\% & 21.58\\
& LMPS-MIRP & 1.77 & 81.26\% & 1.84\% & 18.69\% & 44.20\\
& LMPS-FixedInit & 1.70 & 73.80\% & 2.11\% & 20.00\% & 34.98\\
\multirow{-5}{*}{MSCI} & BAH & 1.37 & 38.58\% & {\color{red} \textbf{1.22}\%} & 12.69\% & 31.71 \\
\hline
\end{tabular}
}
\end{table*}

\begin{figure*}[h!]
    \centering
    \includegraphics[width=\textwidth]{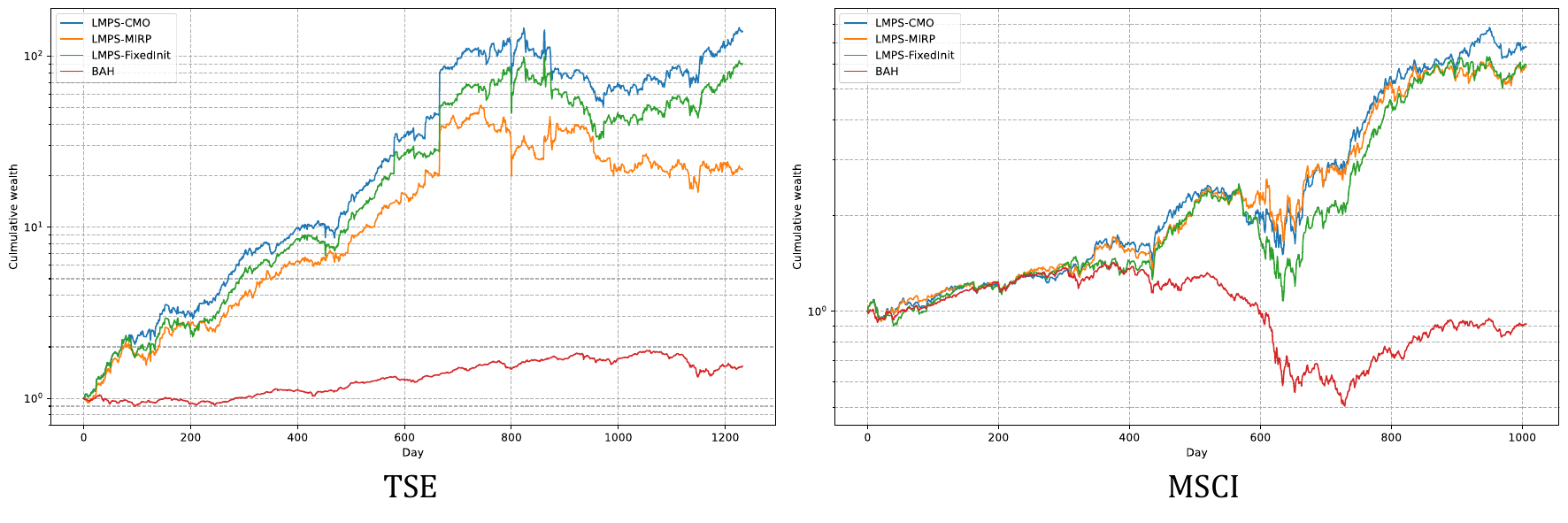}
    \caption{Ablation study of the LMPS-CMO method in TSE and MSCI datasets}
    \label{fig:Ablation studies2}
\end{figure*}

\textbf{Cross-market test:}
Due to the lack of cross-market transferability in the E2E-online method, we only compare LMPS-CMO with LMPS-MIRP, LMPS-FixInit, and BAH in the ablation experiments for cross-market test. The results are displayed in Table \ref{tab:cross_stage_object2} and Figure \ref{fig:Ablation studies2}.

The results indicate that LMPS-CMO outperforms the other three strategies in terms of cumulative returns and Sharpe ratio, effectively mitigating risks. This highlights the robust learning capability of LMPS-CMO. In practical applications, compared to traditional end-to-end models, the proposed model efficiently utilizes data from various markets and historical data.

\begin{table*}[htb!]
\centering
\caption{Ablation study of the LMPS-CMO method in TSE and MSCI datasets}
\label{tab:cross_stage_object2}
\resizebox{\textwidth}{!}{%
\renewcommand{\arraystretch}{1} 
\begin{tabular}{c|cccccccccc}
\hline
Dataset & Strategy & \begin{tabular}[c]{@{}c@{}}Culmulative \\ Wealth\end{tabular} & \begin{tabular}[c]{@{}c@{}}Annualized\\ Return\end{tabular} & \begin{tabular}[c]{@{}c@{}}Annualized\\ Volatility\end{tabular} & \begin{tabular}[c]{@{}c@{}}Maximum\\ Drawdown\end{tabular} & \begin{tabular}[c]{@{}c@{}}Sharpe\\ Ratio\end{tabular} & \\
\hline
\multirow{4}{*}{TSE} &\textbf{LMPS-CMO} & {\color{red} \textbf{149.57}} & {\color{red} \textbf{171.96}\%} & 10.03\% & 65.21\% & {\color{red} \textbf{17.54}} & \\
& LMPS-MIRP & 21.80 & 86.80\% &9.49\% & 68.99\%& 9.15 &  \\
& LMPS-FixedInit & 90.20 & 149.13\% & 10.10\% & 67.20\% & 14.76 &  \\
& BAH & 1.54 & 9.18\% & {\color{red}\bf 1.81\%} & {\color{red} \textbf{30.02}\%} & 5.08 & \\
\hline
\multirow{4}{*}{TSE} &\textbf{LMPS-CMO} & {\color{red} \textbf{6.73}} & {\color{red} \textbf{60.89}\%} & 4.56\% & 39.36\% & {\color{red} \textbf{13.30}} \\
&LMPS-MIRP &  5.86 & 55.11\% & 4.57\% & {\color{red} \textbf{36.40}\%} & 12.05 \\
&LMPS-FixedInit &  5.96 & 55.76\% & 4.70\% & 57.14\% & 11.88 \\
&BAH &  0.91 & -2.22\% & {\color{red} \textbf{3.15}\%} & 64.75\% & -0.71\\
\hline
\end{tabular}}
\end{table*}

\section{Conclusion}
In this paper, we introduce the Meta-LMPS-online algorithm for online portfolio selection problem, applying the meta-learning approach to mixture policies. The core of meta-learning lies in extracting patterns and rules from diverse market data or historical trends, empowering decision-makers to quickly adapt to new investment landscapes. Meanwhile, the mixture policies learning framework enhances stability and versatility by consolidating weighted combinations of multiple candidate polices, thereby refining overall performance. This paper introduces LMPS-SMO strategy for individual markets and LMPS-CMO strategy for multiple markets, effectively addressing the limitations of end-to-end networks in learning from diverse markets. Empirical tests validate
the effectiveness of meta-learning approach and the mixture policies learning framework, highlighting the substantial aid provided by this algorithm in addressing online portfolio selection problems.

\section{Acknowledgments}
This research was supported by National Key R\&D Program of China under No. 2022YFA1004000 and National Natural Science Foundation of China under Grant Number 12371324.



\bibliographystyle{plainnat}
\bibliography{elsarticle}

\section*{Appendix A. Background of Meta-Learning} 

Meta-learning is a newly flourished research direction in the machine learning field. The goal of meta-learning is to learn prior knowledge from many similar tasks such that the learned knowledge can be fast adapted to new unseen tasks. In contrast to traditional  machine learning method involves manually adjusting hyperparameters, meta-learning aims to achieve automatic hyperparameters adjustment. Therefore, meta-learning can be considered as learning beyond the traditional levels of machine learning. The distinction between traditional supervised machine learning and meta-learning is illustrated in Table \ref{the difference}.
\begin{table}[h]
    \centering
    \caption{The difference between machine learning and meta-learning}
    \label{the difference}
    \resizebox{\textwidth}{!}{
    \begin{tabular}{c|c|c}
    \hline
        &Supervised Machine Learning & Meta-Learning\\
        \hline
        Object of fitting&\begin{tabular}[c]{@{}c@{}}Data pairs$\{(x_i,y_i)\}_{i=1}^n$\\ Feature \& Label\end{tabular}&\begin{tabular}[c]{@{}c@{}}Tasks$\{\mathcal{T}_i^{sup},\mathcal{T}_i^{que}\}_{i=1}^m$\\Support set \& Query set \end{tabular}\\
        \hline
        Trained parameters&\begin{tabular}[c]{@{}c@{}}Learning machine parameters \\for predicting data labels\end{tabular}&\begin{tabular}[c]{@{}c@{}}Meta-learning machine parameters for \\predicting machine learning hyperparameters\end{tabular}\\
        \hline
        Evaluation metric&\begin{tabular}[c]{@{}c@{}}Accuracy of predictions\\ on training data labels\end{tabular}&\begin{tabular}[c]{@{}c@{}}Learning performance obtained on the query set\\ by the learning method under the hyperparameters 
        \\assignments obtained from the support set\end{tabular}\\

        \hline
    \end{tabular}}
\end{table}

MAML(Model-Agnostic Meta-Learning) is one of the classic approaches that applies meta-learning to Few-Shot Learning (FSL). Its primary setup involves initializing hyperparameters as the initial values for model training. The objective of the meta-learning machine in MAML is to learn these initial values from a training task set and then directly apply them to prediction tasks. This allows achieving a good training result for the target task after a few iterations. 
The complete training process of MAML is shown in the algorithm \ref{maml}.

\begin{algorithm}[ht]
\caption{MAML}
\label{maml}
\renewcommand{\algorithmicrequire}{\textbf{Input:}}
\renewcommand{\algorithmicensure}{\textbf{Output:}}
\begin{algorithmic}[1]
    \REQUIRE ~~ $p(\mathcal{T})$: distribution over tasks\\
        \REQUIRE ~~ $\alpha,\beta$: step size hyperparameters\\
        \REQUIRE ~~ $K,Q$: the number of samples in the Support Set and Query Set\\
        \REQUIRE ~~ $f_\theta$: the candidate model
        
    \STATE  randomly initialize $\theta$ 
    \WHILE{not done}
            \STATE Sample batch of tasks $\mathcal{T}_i\sim p(\mathcal{T})$
            \FOR{$\mathcal{T}_i$} 
                \STATE Evaluate $\nabla_\theta \mathcal{L}_{\mathcal{T}_i}(f_\theta(\mathcal{T}_i^{sup}))$ with respect to $K$ examples
                \STATE Compute adapted parameters with gradient descent: $\theta_i'=\theta-\alpha \nabla_\theta \mathcal{L}_{\mathcal{T}_i}(f_\theta(\mathcal{T}_i^{sup}))$
            \ENDFOR
            \STATE Update $\theta\leftarrow \theta-\beta\nabla_\theta\sum_{\mathcal{T}_i\sim p(\mathcal{T})} \mathcal{L}_{\mathcal{T}_i}(f_{\theta_i^{'}}(\mathcal{T}_i^{que}))$
        \ENDWHILE
\end{algorithmic}
\end{algorithm}

\end{document}